\documentclass[12pt]{elsarticle}
\usepackage{mathrsfs}
\usepackage{amsmath,amssymb}
\usepackage{lineno,hyperref}
\modulolinenumbers[5]

\newtheorem{thm}{Theorem}
 \newtheorem{lem}[thm]{Lemma}
  \newtheorem{cor}[thm]{Corollary}
    
 \newtheorem{defn}[thm]{Definition}
 
 \newdefinition{rmk}{Remark}
 \newproof{pf}{Proof}
 \newproof{poti}{Proof of Theorem \ref{thm1}}
 \newproof{potii}{Proof of Theorem \ref{thm2}}

\begin{document}

\begin{frontmatter}

\title{$1$-product problems with congruence conditions in nonabelian groups}


\author[mymainaddress]{Kevin Zhao}
\ead{zhkw-hebei@163.com}

\address[mymainaddress]{Department of Mathematics, South China normal university, Guangzhou 510631, China}

\begin{abstract}

  Let $G$ be a finite group and $D_{2n}$ be the dihedral group of $2n$ elements.

  For a positive integer $d$, let $\mathsf{s}_{d\mathbb{N}}(G)$ denote the smallest integer $\ell\in \mathbb{N}_0\cup \{+\infty\}$ such that every sequence $S$ over $G$ of length $|S|\geq \ell$ has a nonempty $1$-product subsequence $T$ with $|T|\equiv 0$ (mod $d$).
  In this paper, we mainly study the problem for dihedral groups $D_{2n}$ and determine their exact values:
  $\mathsf{s}_{d\mathbb{N}}(D_{2n})=2d+\lfloor log_2n\rfloor$, if $d$ is odd with $n|d$;
$\mathsf{s}_{d\mathbb{N}}(D_{2n})=nd+1$, if $gcd(n,d)=1$.
Furthermore, we also analysis the problem for metacyclic groups $C_p\ltimes_s C_q$ and obtain a result:
$\mathsf{s}_{kp\mathbb{N}}(C_p\ltimes_s C_q)=lcm(kp,q)+p-2+gcd(kp,q)$,
where $p\geq 3$ and $p|q-1$.

\end{abstract}


\end{frontmatter}

\section{Introduction}

Let $G$ be a finite group and $D_{2n}$ be the dihedral group of $2n$ elements.
For a finite abelian group $G$,
it is well known that $|G|=1$ or $G = C_{n_1} \oplus C_{n_2}\oplus \ldots
\oplus C_{n_k}$ with $1 < n_1|n_2| \ldots |n_k$,
where $r(G)=k$ is the rank of $G$ and the exponent
$\exp(G)$ of $G$ is $n_k$.
Set $$\mathsf{D}^*(G):=1+\sum _{i=1} ^{r} (n_i-1).$$
For a finite group $G$, we define $\exp(G)$ to be the maximum value of $ord(g)$ over all $g\in G$.
Obviously, the definition covers the one of exponent of a finite abelian group.
The Davenport constant $\mathsf{D}(G)$ is the minimal integer $\ell \in \mathbb{N}_0\cup \{+\infty\}$
such that every sequence $S$ over $G$ of length $|S| \geq \ell$
has a nonempty $1$-product subsequence.
It is easy to see that $\mathsf{D}(G)\geq \mathsf{D}^*(G)$ for any finite abelian group $G$.
Let $\mathsf{s}_{k\exp(G)}(G)$ denote the smallest integer $\ell\in \mathbb{N}_0\cup \{+\infty\}$ such that every
sequence $S$ over $G$ of length $|S|\geq \ell$ has a nonempty
$1$-product subsequence $T$ of length $|T|=k\exp(G)$, where $k$ is a positive integer.
In this paper, we investigate a following generalization of $\mathsf{D}(G)$ and $\mathsf{s}_{k\exp(G)}(G)$.

\begin{defn}
Let $G$ be a finite group.
For a positive integer $d$, let $\mathsf{s}_{d\mathbb{N}}(G)$ denote the smallest integer $\ell\in \mathbb{N}_0\cup \{+\infty\}$ such that every sequence $S$ over $G$ of length $|S|\geq \ell$ has a nonempty $1$-product subsequence $T$ with $|T|\equiv 0$ (mod $d$).
\end{defn}

The invariant $\mathsf{s}_{d\mathbb{N}}(G)$ was introduced by A. Geroldinger etc. \cite{[GGS]}.
It is trivial to see that $\mathsf{s}_{d\mathbb{N}}(G)=\mathsf{D}(G)$, if $d=1$;
$\mathsf{s}_{d\mathbb{N}}(G)\leq \mathsf{s}_{k\exp(G)}(G)$, if $d=k\exp(G)$.
In addition, it is easy to see that
$$S=1^{[d-1]}S_1\in \mathscr{F}(G)$$
is $d\mathbb{N}$-$1$-product free, where $S_1$ is a $1$-product free subsequence of length $D(G)-1$.
Hence,
\begin{equation}
\label{low}
   \begin{aligned}
   &\mathsf{s}_{d\mathbb{N}}(G)\geq D(G)+d-1
   \end{aligned}
  \end{equation}
for any positive integer $d$.
In particular, for $d=k\exp(G)$ with $k\in \mathbb{N}$ we have
\begin{equation}
\label{upper}
   \begin{aligned}
   &k\exp(G)+D(G)-1\leq \mathsf{s}_{k\exp(G)\mathbb{N}}(G)\leq \mathsf{s}_{k\exp(G)}(G).
   \end{aligned}
  \end{equation}
For general, the problem of determining $\mathsf{s}_{d\mathbb{N}}(G)$ is not at all trivial.
Recently, A. Geroldinger etc. \cite{[GGS]} determined the exact values of $\mathsf{s}_{d\mathbb{N}}(G)$ for all $d\geq 1$ when $G$ is a finite abelian group with rank at most two.
That is, $$\mathsf{s}_{d\mathbb{N}}(G)= D^*(G\oplus C_d) = lcm(n,d)+gcd(n,lcm(m,d))+gcd(m,d)-2$$ for $G\cong C_m\oplus C_n$ with $1\leq m | n$, and
$$\mathsf{s}_{d\mathbb{N}}(G)= D^*(G\oplus C_d) = lcm(n,d)+gcd(n,d)-1$$ for $G\cong C_n$.
From the above results it follows that $\mathsf{s}_{d\mathbb{N}}(G)$ can take the lower bound of (\ref{low}).
Furthermore, they also obtained precise values in the case of $p$-groups under mild conditions on $d$.
They found that most of these values can also take the lower bound of (\ref{low}).

We are interested in the case $d=k\exp(G)$ with $k\in \mathbb{N}$.
From the above we can see that $\mathsf{s}_{kn\mathbb{N}}(C_n)= (k+1)n-1=kn+n-1=\mathsf{s}_{kn}(C_n)$.
In \cite{[GHPS]}, the authors proved the following theorem:

\begin{thm} [\cite{[GHPS]},  Theorem 1.1] \label{impor}
Let $H$ be an arbitrary finite abelian group with $exp(H) = m \geq 2$,  and let $G = C_{mn} \bigoplus H$.  If $n \geq 2m|H| + 2|H|$, then $\mathsf{s}_{kmn}(G) = kmn+ D(G)-1$ for all positive integers $k\geq 2$.
\end{thm}
Combining Theorem \ref{impor} with (\ref{upper}) yields the following results:

\begin{cor} \label{dNcong}
Let $H$ be an arbitrary finite abelian group with $exp(H) = m \geq 2$,  and let $G = C_{mn} \bigoplus H$.  If $n \geq 2m|H| + 2|H|$, then $\mathsf{s}_{kmn\mathbb{N}}(G) = kmn+ D(G)-1$ for all positive integers $k\geq 2$.
\end{cor}
Hence, $\mathsf{s}_{d\mathbb{N}}(G)$ can take the upper bound of (\ref{upper}).

In this paper, we mainly study the invariant $\mathsf{s}_{d\mathbb{N}}(D_{2n})$ with $d$ odd and $n|d$.
Since $d$ is odd, then it is easy to see that $s_{d}(D_{2n})=+\infty$,
since $x^{[N]}$ is $d$-$1$-product free for any $N\in \mathbb{N}$.
Thus from Theorem \ref{main} it follows that
$$d+D(D_{2n})-1=d+n< \mathsf{s}_{d\mathbb{N}}(D_{2n})< +\infty=\mathsf{s}_{k\exp(G)}(G),$$
that is, $\mathsf{s}_{d\mathbb{N}}(G)$ can take the middle values of (\ref{upper}).

We use the following generators and relations for the dihedral group of order $2n$:
$$D_{2n}=\langle x,y:x^2=y^n=1, xy=y^{-1}x\rangle.$$
Following the notations of \cite{[B]} and \cite{[GL]},
let $H$ be the cyclic subgroup of $D_{2n}$ generated by $y$,
and let $N=D_{2n}\setminus H.$
It is easy to see that $N=D_{2n}\setminus H=x\cdot H.$
For any sequence $S\in\mathscr {F}(D_{2n})$, we denote $S\cap H$ and $S\cap N$ by $S_H$ and $S_N$, respectively.

Our main result is the following:
\begin{thm} \label{main}
Let $n$, $d$ be positive integers, $n\geq 3$.
If $d$ is odd and $n|d$, then
$$\mathsf{s}_{d\mathbb{N}}(D_{2n})=2d+\lfloor log_2n\rfloor.$$

Furthermore, if $gcd(n,d)=1$, then $\mathsf{s}_{d\mathbb{N}}(D_{2n})=nd+1.$
\end{thm}

We also analysis the problem for metacyclic groups $G_{pq}=C_p\ltimes_s C_q$ and obtain a result:
\begin{thm} \label{thm2}
Let $G_{pq}=\langle x,y: x^p=y^q = 1,yx = xy^s,ord_q(s) =p,$ and $p,q$ primes$\rangle$ be an nonabelian group.
If $p\geq 3$ and $p|q-1$, then
$$
\mathsf{s}_{kp\mathbb{N}}(G_{pq})=lcm(kp,q)+p-2+gcd(kp,q)
$$

\end{thm}

\section{Preliminaries}

Now we present some key concepts.
Let $\mathbb{N}$ denote the set of positive integers and $\mathbb{N}_0=\mathbb{N} \cup \{0\}$.
For real numbers $a,b\in \mathbb{R}$, we set $[a,b] = \{x\in \mathbb{Z} : a\leq x\leq b\}$.
For any $x\in \mathbb{R}$, denote $\lfloor x\rfloor$ by the smallest integer $t\leq x$.
For positive integers $m,\ n$, denote by $gcd(m,\ n)$ and $lcm(m,\ n)$ the greatest common divisor and
the least common multiple of $m,\ n$ respectively.
For positive integers $n$ and $g$ with $(n,g)=1$,
let $ord_n(g)$ denote the minimal positive integer $\ell$ such that $g^{\ell}\equiv 1$ (mod $n$).

Let $G$ be a finite group with multiply operation $\cdot$ and unit element $1$
and $C_n$ be the cyclic group.
If $H$ is a subgroup of $G$ and $H\neq G$, then we call $H$ a proper subgroup of $G$.
Let $G^{\bullet}=G\setminus \{1\}$.
Let $\mathscr {F}(G)$ be the free abelian monoid,
multiplicatively written, with basis $G$.
The elements of $\mathscr {F}(G)$ are called sequences over $G$.
Let
$$S=g_1 \cdot\ldots\cdot g_{\ell}= {\prod_{g\in G}}^{\bullet}g^{[v_g(S)]}\in \mathscr{F}(G)$$
with $v_g(S) \in \mathbb{N}_0$ for all $g \in G$.
We call $v_g(S)$ the multiplicity of $g$ in $S$,  and if $v_g(S) > 0$ we say that $S$ contains $g$.
We call $S$ a \textit{$\pm$-$1$-product sequence} (\textit{$1$-product sequence})
if $1=\prod_{i=1}^{\ell}g_{\tau(i)}^{\varepsilon_i}$ ($1=\prod_{i=1}^{\ell}g_{\tau(i)}$) holds for some permutations $\tau$ of $\{1,\ldots ,l\}$, where $\varepsilon_i\in \{1,-1\}$ for each $i\in [1,{\ell}]$.
We denote by $\Pi_{\pm}(S)$ ($\Pi(S)$) the set of all $\pm$-products (products) $\prod_{i=1}^{\ell}g_{\tau(i)}^{\varepsilon_i}$ ($\prod_{i=1}^{\ell}g_{\tau(i)}$), where $\tau$ runs over all permutation of $[1,\ell]$,
i.e., $\Pi_{\pm}(S)=\{\prod_{i=1}^{\ell}g_{\tau(i)}^{\varepsilon_i} : \tau$ $is\ a\ permutation\ of\ [1,\ell],\ \varepsilon_{i}\in \{1,-1\}\}$ ($\Pi(S)=\{\prod_{i=1}^{\ell}g_{\tau(i)} : \tau$ is a permutation of $[1,\ell]\}$).
Note that if $(G,+,0)$ is an additive finite abelian group, then we can replace $\pm$-$1$-product by $\pm$-zero-sum.
Definite $\pi_{\pm}(S)=g_1^{\varepsilon_1}\cdot g_2^{\varepsilon_2}\cdot\ldots\cdot g_{\ell}^{\varepsilon_{\ell}}$ ($\pi(S)=g_1g_2\cdot \ldots\cdot g_{\ell}$) to be the specific $\pm$-product (product) of $S$ obtained by multiplying all elements in the order they appear in $S$.
Denote $\mathsf{D}_{\pm}(G)$ by the minimal integer $\ell \in \mathbb{N}$
such that every sequence $S$ over $G$ of length $|S| \geq \ell$
has a nonempty $\pm$-$1$-product subsequence.

If for all $g \in G$ we have $v_g(S) = 0$, then we call $S$ the \textit{empty sequence}.
Apparently, a squarefree sequence over $G$ can be treated as a subset of $G$.
A sequence $S$ is called \textit{squarefree} if $v_g(S) \leq 1$ for all $g \in G$.
A sequence $S_1 \in \mathscr {F}(G)$ is called a subsequence of $S$ if $v_g(S_1) \leq v_g(S)$ for all $g \in G$,
and denoted by $S_1 \mid S$.
It is called a \textit{proper subsequence} of $S$ if it is a subsequence with $1 \neq S_1 \neq S$.
Let $S_1, S_2 \in \mathscr {F}(G)$,  we set
\begin{equation*}
S_1\cdot S_2={\prod_{g\in G}}^{\bullet}g^{[v_g(S_1)+v_g(S_2)]}\in \mathscr {F}(G).
\end{equation*}
\begin{equation*}
S_1\cap S_2={\prod_{g\in G}}^\bullet g^{[min\{v_g(S_1),v_g(S_2)\}]}\in \mathscr {F}(G).
\end{equation*}
\begin{equation*}
S_1S_2^{-1}={\prod_{g|S_1}}^\bullet g^{[max\{0,v_g(S_1)-v_g(S_2)\}]}\in \mathscr {F}(G).
\end{equation*}
In particular, if $S_2 \mid S_1$,  we set
\begin{equation*}
S_1S_{2}^{-1}={\prod_{g\in G}}^{\bullet}g^{[v_g(S_1)-v_g(S_2)]}\in \mathscr {F}(G).
\end{equation*}
For the sequence $S$,
we list the following definitions:
\begin{align*}
&|S|=\ell=\sum_{g\in G}v_g(S)\in \mathbb{N}_0,\mbox{ the length of $S$};  \\
&h(S) = \max\{v_g (S):\; g\in G\} \in [0, |S|],\mbox{ the maximum multiplicity of $S$};  \\
&supp(S) = \{g \in G :\; v_g(S) > 0\} \subseteq G,\mbox{ the support of $S$};  \\
&\Sigma(S) = \{\prod\limits \limits _{i\in I}g_{\tau(i)} : I\subset [1, \ell]\ with\ 1 \leq |I| \leq \ell\ and \ \tau \ is\ a\ permutation\ of\ I\}, \\
&\mbox{ the set of all subproducts of $S$}; \\
&\Sigma_{\pm}(S) = \{\prod\limits \limits _{i\in I}g_{\tau(i)}^{\varepsilon_i} : \varepsilon_i\in \{1,-1\}, \ I\subset [1, \ell]\ with\ 1 \leq |I| \leq \ell\ and \ \tau \ is\ a\ permutation\ of\ I\}, \\
&\mbox{ the set of all $\pm$-subproducts of $S$}; \\
&\Sigma_k(S) = \{\prod_{i\in I}g_{\tau(i)} : I\subseteq [1, \ell]\ with\ |I|=k \ and \ \tau \ is\ a\ permutation\ of\ I\}, \\
&\mbox{ the set of $k$-term subproducts of $S$}, \\
&\Sigma_{d\mathbb{N}}(S) = \{\prod_{i\in I}g_{\tau(i)} :\; I\subseteq [1, \ell]\ with\ d||I|\ and \ \tau \ is\ a\ permutation\ on\ I\}. \\
\end{align*}
In addition, we write $${\sum} _{\leq k}(S) = \bigcup \limits _{j \in [1, k]}{\sum}_j(S),
{\sum} _{\geq k}(S) = \bigcup \limits_{j \geq k}{\sum}_j(S),$$
and
$${\sum}_E(S)=\{\sum(T): T|S \ and \ 2||T|\},
{\sum}_O(S)=\{\sum(T):T|S \ and \ 2\nmid |T|\}.$$
The sequence $S$ is called
\begin{itemize}
  \item \textit{$1$-product free} if $1\not\in \Sigma(S)$,
  \item \textit{$n$-$1$-product free} if $1\not\in \Sigma_n(S)$,
  \item \textit{$d\mathbb{N}$-$1$-product free} if $1\not\in \Sigma_{d\mathbb{N}}(S)$,
  \item \textit{a minimal $\pm$-$1$-product ($1$-product) sequence} if $S$ is nonempty, $1\in \sum_{\pm}(S)$ ($1\in \sum(S)$), and every $S'\mid S$ with $1\leq |S'|< |S|$ is $\pm$-$1$-product ($1$-product) free.
\end{itemize}

In the process of the proof, we will use some basic results in the zero-sum theory.

\begin{lem} \label{basic}
Let $m$, $n$, $d$ be positive integers, $m|n$, $n\geq 2$ and let $G$ be a finite abelian group,
$C_n$ be a cyclic group of order $n$, $D_{2n}$ be the dihedral group. Then

\begin{description}
  \item[(a)] $D(C_n)=n$. If $S\in \mathscr {F}(C_n)$ is a minimal $1$-product sequence of length $n$, then $S=g^{[n]}$ for some $g\in G$ with $ord(g)=n$.(\cite{[GH]}, Theorem 5.1.10.1)
  \item[(b)] $s_{dn}(C_n)=(d+1)n-1$. (\cite{[GH]}, Corollary 5.7.5)
  \item[(c)] $D(C_{m}\oplus C_{n})=D^*(C_{m}\oplus C_{n})=m+n-1$. (\cite{[GH]}, Theorem 5.8.3)
  \item[(d)] $C_n\oplus C_d\cong C_{gcd(n,d)}\oplus C_{lcm(n,d)}$. (\cite{[GGS]}, Lemma 3.2)
  \item[(e)] If $D^*(G\oplus C_d)=D(G\oplus C_d)$, then $s_{d\mathbb{N}}(G)=D^*(G\oplus C_d)$.
  (\cite{[GGS]}, Proposition 3.3)
  \item[(f)] $D(D_{2n})=n+1$. (\cite{[B]}, Lemma 4)
  \item[(g)] $s_{2n}(D_{2n})=3n$. (\cite{[B]}, Theorem 8)
\end{description}
\end{lem}

\begin{lem} [\cite{[MBR]}, Theorem 1.3] \label{inverse-dihedral}
Let $S$ be a sequence of length $n$ in the dihedral group $D_{2n}$, where $n\geq 3$.
\begin{description}
  \item[(1)] If $n\geq 4$, then the following statements are equivalent:

$(i)$ $S$ is $1$-product free;

$(ii)$ $S = (y^t)^{[n-1]}, xy^s$
 with $1\leq t\leq n-1$, $gcd(t,n) = 1$ and $0\leq s\leq n-1$.
  \item[(2)] If $n\geq 3$, then $S$ is $1$-product free if and only if
either $S = (y^t,y^t,xy^{\nu})$ with $t\in \{2,3\} \ and \ \nu \in \{0,1,2\} \ or \ S = (x,xy,xy^2).$
\end{description}
\end{lem}

By the above lemmas, we can easily show following results:
\begin{lem} \label{lem40}
\begin{description}
  \item[(i)] Let $m$ be a positive integer, and let $G$ be a finite abelian group.
  Suppose that $S\in \mathscr {F}(G)$ is a $1$-product free sequence of length $m$.
  Then $|\sum(S)|\geq |S|$. In particular, if $|\sum(S)|=|S|$, then $S=g^{[m]}$ for some $g\in G$ with $ord(g)>m$.
  \item[(ii)] Let $n\geq 3$ be a positive integer and let $D_{2n}$ be the dihedral group. Then
  $$s_{\leq n}(D_{2n})=n+1.$$
\end{description}

\end{lem}

\pf
$(i)$ Set $$S=g_1\cdot \ldots\cdot g_m\in \mathscr {F}(G)$$
is a $1$-product free sequence.
It is easy to see that $g_1,g_1g_2,g_1g_2g_3,\ldots ,g_1g_2\cdot\ldots\cdot g_m$ are distinct,
which implies that $|\sum(S)|\geq |S|$.
Suppose that $|supp(S)|\geq 2$ and let $g_1\neq g_2$.
Since $S$ is $1$-product free, we have that $g_1,g_2,g_1g_2,$$g_1g_2g_3,\ldots ,g_1g_2\cdot\ldots\cdot g_m$ are all distinct and contained in $\sum(S)$.
This is in contradiction to $|\sum(S)|=|S|$.
Hence, $|supp(S)|=1$ and then $S=g^{[m]}$ for some $g\in G$ with $ord(g)>m$.

$(ii)$ Let $D_{2n}=\langle x,y:x^2=y^n=1, xy=y^{-1}x\rangle$.
It is easy to see that
$$S'=xy^{[n-1]}$$
is a $1$-product free sequence of length $n$.
Thus it suffices to prove that $s_{\leq n}(D_{2n})\leq n+1$.

Fix a sequence $S$ of length $n+1$.
If $|S_H|\geq n$, then Lemma \ref{basic} $(a)$ completes the proof.
If $|S_H|\leq n-1$, then $|S_N|=|S|-|S_H|\geq 2$.
Take a subsequence $T|S$ of length $n$ satisfying $|T_N|\geq 2$.
From Lemma \ref{inverse-dihedral} it follows that
if $T\in \mathscr {F}(D_{2n})$ is a $1$-product free sequence of length $n$, then $|T_N|=1$.
Therefore, there exists a $1$-product subsequence $T_1|T$ of length $|T_1|\leq n$,
and we complete the proof.

\qed

The following lemmas and corollaries are repeatedly used.

\begin{lem}[\cite{[GH]}] \label{maxlength}
Let $G$ be a finite abelian group of order $n$, and let $S$ be a sequence of $n$ elements in $G$. Let $k =max\{v_g(S):g\in G\}$ be the maximal value of repetition of an element occurring in $S$. Then $1\in \sum_{\leq k}(S)$.
\end{lem}

\begin{lem}[\cite{[AD]}, Lemma 2.1] \label{cyweight}
Let $n\in \mathbb{N}$ and $(y_1, \ldots ,y_s)$ be a sequence of integers with $s>log_2n$. Then there exists a nonempty $J \subseteq \{1,2,3, \ldots ,s\}$ and $\varepsilon_j\in \{\pm 1\}$ for each $j\in J$ such that

$$
\sum_{j\in J}\varepsilon_jy_j\equiv 0 \ (mod \ n)
$$

\end{lem}

\begin{cor} \label{dihedral}
Suppose that $T_1\in \mathscr {F}(N)$, $T_2\in \mathscr {F}(H)$ are two sequences satisfying that $\lfloor\frac{|T_1|}{2}\rfloor$, $\lfloor\frac{|T_2|}{2}\rfloor >\lfloor log_2n\rfloor$.
Then
\begin{description}
  \item[(i)] there exists a $1$-product subsequence $T'_1$ in $T_1$ satisfying that $2||T'_1|$ and $|T'_1|\leq 2\lfloor log_2n\rfloor+2$;
  \item[(ii)] there exist two disjointing subsequences $W_1$, $W_2$ in $T_2$ satisfying that $1\leq |W_1|=|W_2|\leq \lfloor log_2n\rfloor+1$ and $\pi(W_1)=\pi(W_2)$.
\end{description}
\end{cor}

\pf The proof of the Corollary is from \textit{Proof of Theorem $1.3$} of \cite{[MBR]}.
For the convenience, we exhibit it here.

$(i)$. Set $T_1={\prod^{\bullet}}_{i\in [1,k]}xy^{\alpha_i}.$
Thus $k=|T_1|$ and $\lfloor\frac{k}{2}\rfloor >\lfloor log_2n\rfloor.$
From Lemma \ref{cyweight} it follows that there exists a linear combination of a subset $T''_1$ of
$$
\{(\alpha_1-\alpha_2), (\alpha_3-\alpha_4), \ldots , (\alpha_{2\lfloor\frac{k}{2}\rfloor-1}-\alpha_{2\lfloor\frac{k}{2}\rfloor})\}
$$
with coefficients $\pm1$ summing $0$ satisfying that $v=|T''_1|\leq\lfloor log_2n\rfloor+1$.
Suppose without loss of generality that, in this combination
$$
(\alpha_1-\alpha_2), (\alpha_3-\alpha_4), \ldots , (\alpha_{2u-1}-\alpha_{2u})
$$
appear with signal $-1$ and
$$
(\alpha_{2u+1}-\alpha_{2u+2}), (\alpha_{2u+3}-\alpha_{2u+4}), \ldots , (\alpha_{2v-1}-\alpha_{2v}),
$$
appear with signal $+1$. Then
\begin{align*}
(xy^{\alpha_1}\cdot xy^{\alpha_2})\cdot \ldots\cdot (xy^{\alpha_{2u-1}}\cdot xy^{\alpha_{2u}})\cdot (xy^{\alpha_{2u+1}}\cdot  xy^{\alpha_{2u+2}})\cdot \ldots\cdot (xy^{\alpha_{2v-1}}\cdot xy^{\alpha_{2v}})= \\
y^{\alpha_2-\alpha_1}\cdot \ldots\cdot y^{\alpha_{2u}-\alpha_{2u-1}}\cdot y^{\alpha_{2u+2}-\alpha_{2u+1}}\cdot \ldots\cdot y^{\alpha_{2v}-\alpha_{2v-1}}= 1
\end{align*}
and we find a $1$-product subsequence satisfying the conditions of $(i)$.

$(ii)$. By imitating the proof of $(i)$, one can obtain $(ii)$.
One can also find the proof in Lemma $8$ of \cite{[GL]}.

\qed

\begin{lem} \label{dihedral-length}
Let $n\geq 2$ be a positive integer. Then
$$D_{\pm}(D_{2n})\leq 2\lfloor log_2n\rfloor+2.$$
Furthermore, if $S\in \mathscr {F}(D_{2n})$ is a $\pm$-$1$-product sequence with $|S_N|>0$,
then $S$ is $1$-product.
\end{lem}

\pf
Let $D_{2n}=\langle x,y:x^2=y^n=1, xy=y^{-1}x\rangle$, $H=\langle y\rangle$ and $N=x\cdot H$.
Suppose that
$$S=xy^{\alpha_1}, \ldots ,xy^{\alpha_t},y^{\beta_1}, \ldots ,y^{\beta_k}\in \mathscr {F}(D_{2n})$$
is a sequence with $|S|=t+k\geq 2\lfloor log_2n\rfloor+2.$
Set $$S_0=(\alpha_1-\alpha_2), (\alpha_3-\alpha_4), \ldots , (\alpha_{2\lfloor\frac{t}{2}\rfloor-1}-\alpha_{2\lfloor\frac{t}{2}\rfloor}),\beta_1, \ldots ,\beta_k.$$
Thus, $|S_0|=\lfloor\frac{t}{2}\rfloor+k\geq \lfloor log_2n\rfloor+1$.
By imitating the proof of Corollary \ref{dihedral},
we can find a $\pm$-$1$-product subsequence in $S$.
The first assertion is c1ear.

It will not be too confusing to suppose that
$$S=xy^{\alpha_1}, \ldots ,xy^{\alpha_t},y^{\beta_1}, \ldots ,y^{\beta_k}\in \mathscr {F}(D_{2n})$$
is a $\pm$-$1$-product sequence with $|S_N|>0$.
In addition, it is no loss of generality to assume that
$$\prod_{i=1}^{t}xy^{\alpha_i}\prod_{j\in I}y^{\beta_j}\prod_{\ell\in \overline{I}}y^{-\beta_{\ell}}=1$$
where $I\subseteq [1,k]$ and $\overline{I}=[1,k]\setminus I$.
Thus $$\prod_{i=1}^{t-1}xy^{\alpha_i}\prod_{\ell\in \overline{I}}y^{\beta_{\ell}}xy^{\alpha_t}\prod_{j\in I}y^{\beta_j}=1.$$
The proof of the second assertion is complete.

\qed

\begin{lem}[\cite{[GHQQZ]}] \label{cyclicadd}
Let $C_n$ be a cyclic group of order $n\in \mathbb{N}$
and let $S\in \mathscr {F}(C_{n})$ be a sequence of length $|S|\geq n$.
If $|S_H|\leq |H|-1$ for any proper subgroup $H$ of $C_n$, then

$$
\sum(S)=C_n.
$$

\end{lem}

\begin{lem}[\cite{[QH]}] \label{cy-additive-inv}
Let $G$ be a cyclic group of order $n$, and let $S$ be a sequence of length $n-1$ in $G$.
If $|S_H|\leq |H|-1$ for any proper subgroup $H$ of $G$, and $\sum(S)\neq G$,
then $S=g^{[n-1]}$ where $g$ is a generator of $G$.
\end{lem}

By Lemma \ref{cyclicadd} and Lemma \ref{cy-additive-inv},
we can prove the following lemma which is crucial in this paper.
\begin{lem} \label{ODD}
Let $n$ be a positive integer with $2\nmid n$, and let $C_n$ be a cyclic group of order $n$.
If $S\in \mathscr {F}(C_{n})$ is a sequence of length $|S|\geq n$,
then $S$ has a $\pm$-$1$-product subsequence $T$ with $2\nmid |T|$.

Furthermore, if $S$ is a minimal $\pm$-$1$-product subsequence of length $n$,
then $S=g^{[n]}$, where $g$ is a generator of $C_n$.
\end{lem}

\pf
If $1$ is contained in $S$, then $T=1$ is $\pm$-$1$-product with $2\nmid |T|$.
Now suppose that $1$ is not contained in $S$.

\textit{Case $1$}: For any proper subgroup $H$ of $C_n$, $|S_H|\leq |H|-1$.

Then by $|S|\geq n$, it is easy to see that $S$ contains an element $g_0$ with $ord(g_0)=n$.
Set $$S_1=Sg_0^{-1}.$$
It follows that $|S_1|\geq n-1$.
If $\sum_E(S)\cap \sum_O(S)\neq \phi$,
then $S$ has two nonempty subsequences $T_1$ and $T_2$ satisfying that $2||T_1|$, $2\nmid |T_2|$ and $\pi(T_1)=\pi(T_2)$.
Let $T'_1=T_1T^{-1}_2$ and $T'_2=T_2T^{-1}_1$.
Thus one can obtain that $\pi(T'_1)=\pi(T'_2)$ and $T'_1T'_2$ is a subsequence of $S$ satisfying that
$$|T'_1T'_2|=|T_1|+|T_2|-2|T_1\cap T_2|$$ is odd.
It follows that $T'_1T'_2$ is $\pm$-$1$-product with $2\nmid |T'_1T'_2|$.
If $\sum(S_1)\neq C_n$, then by Lemma \ref{cyclicadd} we have $|S_1|\leq n-1$.
From the choice of $S_1$ it follows that $|S_1|=n-1$.
Combining the above assumptions with Lemma \ref{cy-additive-inv} yields that $S_1=Sg_0^{-1}=g^{[n-1]}$, where $g$ is a generator of $C_n$.
Set $g_0=g^x$ with $x\in [1,n-1]$.
Let $S'=g_0g^{[n-x]}$ if $2\nmid x$; otherwise let $S'=g_0g^{[x]}$.
For $2\nmid n$, we have that $S'$ is a $\pm$-$1$-product subsequence of odd length.

Now suppose that $\sum_E(S)\cap \sum_O(S)=\phi$ and $\sum(S_1)=C_n$.
Thus $${\sum}_E(S_1)\cap {\sum}_O(S_1)=\phi$$ and
$$C_n=\sum(S)={\sum}_E(S)\cup {\sum}_O(S)=\sum(S_1)={\sum}_E(S_1)\cup {\sum}_O(S_1).$$
It follows that $$|{\sum}_E(S_1)|+|{\sum}_O(S_1)|=|{\sum}_E(S)|+|{\sum}_O(S)|=n.$$
In addition, we can suppose that $g_0\in \sum_O(S)$,
since otherwise $g_0\in \sum_E(S)$, and it follows that
there exists a subsequence $S'|S$ of even length such that $S'g$ is a $\pm$-$1$-product subsequence of odd length.
Since $\sum_E(S_1)\subseteq \sum_E(S)$ and $\sum_O(S_1)\subseteq \sum_O(S)$,
we must have that $${\sum}_E(S_1)={\sum}_E(S) \ and \ {\sum}_O(S_1)={\sum}_O(S).$$
Since ${\sum}_O(S)={\sum}_O(S_1)\cup (g_0\cdot {\sum}_E(S_1))\cup \{g_0\}={\sum}_O(S_1)\cup (g_0\cdot {\sum}_E(S_1))$
and ${\sum}_E(S)={\sum}_E(S_1)\cup (g_0\cdot {\sum}_O(S_1)),$
we have $$g_0\cdot {\sum}_E(S_1)\subseteq {\sum}_O(S_1) \ and \ g_0\cdot {\sum}_O(S_1)\subseteq {\sum}_E(S_1).$$
Thus $|\sum_E(S_1)|=|g_0\cdot \sum_E(S_1)|\leq |\sum_O(S_1)|$
and $|\sum_O(S_1)|=|g_0\cdot \sum_O(S_1)|\leq |\sum_E(S_1)|$.
It follows that $|\sum_E(S_1)|=|\sum_O(S_1)|,$
which implies that $n=2|\sum_E(S_1)|$ is even.
This is in contradiction to $2\nmid n$.

\textit{Case $2$}: There exists a proper subgroup $N$ in $C_n$ such that $|S_{N}|\geq |N|$.

Then the number of divisors of $|N|$ is less than that of $n$.
The following proof is by induction on the number $k$ of divisors of $n$.

If $k=2$, then $n$ is a prime.
It follows that the unique proper subgroup $H$ of $C_n$ is $\{1\}$.
For $1\nmid S$, we have $|S_H|\leq |H|-1$.
By Case $1$, the proof is complete.
Now suppose that the lemma has been established for integers less than $k$.
By the induction hypothesis, there exists a subsequence $T|S_{N}$ such that $T$ is a $\pm$-$1$-product subsequence of odd length.
We complete the proof of the first assertion.

Let $S\in \mathscr {F}(C_{n})$ be a minimal $\pm$-$1$-product subsequence of length $n$.
Suppose that $S$ has a decomposition:
$$S=T_1T_2$$ with $1\leq |T_1|\leq |T_2|$ and $\pi(T_1)=\pi(T_2)$.
For $n$ odd, we must have $|T_1|\neq |T_2|$, i.e., $1\leq |T_1|< |T_2|$.
By the minimality of $S$, we have that $T_1$ and $T_2$ are $1$-product free and $\sum(T_1)\cap \sum(T_2)=\pi(T_1)=\pi(T_2)$.
It follows that $1\not\in\sum(T_1)\cup \sum(T_2)\subseteq C_n$, $|\sum(T_1)|\geq |T_1|$ and $|\sum(T_2)|\geq |T_2|$.
Combining the above results yields that
$n-1\geq|\sum(T_1)\cup \sum(T_2)|=|\sum(T_1)|+ |\sum(T_2)|-|\sum(T_1)\cap \sum(T_2)|\geq |T_1|+|T_2|-1=|S|-1=n-1$,
that is, $|\sum(T_1)|=|T_1|$ and $|\sum(T_2)|=|T_2|$.
By Lemma \ref{lem40} $(i)$, we have that $T_1=g_1^{[|T_1|]}$ and $T_2=g_2^{[n-|T_1|]}$
for some $g_1,\ g_2\in C_n$ with $ord(g_1)>|T_1|$ and $ord(g_2)>n-|T_1|$.
For $\pi(T_1)=\pi(T_2)$,
we have $g_1^{|T_1|}=g_2^{n-|T_1|}$, i.e., $g_1^{|T_1|}\cdot g_2^{|T_1|}=1$.
Combining this with $1\leq |T_1|< |T_2|=n-|T_1|$ yields that $g_1^{[|T_1|]} g_2^{[|T_1|]}$
is a proper $\pm$-$1$-product subsequence of $S$, a contradiction.
Hence, $S$ is a minimal $1$-product subsequence of length $n$.
By Lemma \ref{basic} $(a)$, we complete the proof.

\qed

\section{The proof of Theorem \ref{main}}

\begin{lem}[\cite{[GL]}, Lemma 7] \label{lem12}
Let $G$ be a finite abelian group of order $n$ and let $r \geq 2$ be an integer.
Suppose that $S$ is a sequence of $n+r-2$ elements in $G$.
If $1\notin \sum_n(S)$, then $|\sum_{n-2}(S)|=|\sum_r(S)|\geq r-1$.
\end{lem}

\begin{lem} [\cite{[N]}, Lemma 2.2] \label{lem11}
Let $A$, $B$ be two subsets of a finite group $G$. If $|A|+|B| > |G|$, then $A+ B = G$, where $A+ B =\{ab:a \in A,b\in B\}$.
\end{lem}

\begin{thm} \label{n-dihedral}
Let $n$ be a positive integer, $n\geq 3$.
If $n$ is odd, then
$$s_{n\mathbb{N}}(D_{2n})=2n+\lfloor log_2n\rfloor.$$
\end{thm}

\pf
Let $D_{2n}=\langle x,y:x^2=y^n=1, xy=y^{-1}x\rangle$, $H=\langle y\rangle$ and $N=x\cdot H$.
In addition, for any sequence $S\in\mathscr {F}(D_{2n})$, we denote $S\cap H$ and $S\cap N$ by $S_H$ and $S_N$, respectively.
Set $W=x^{[2n-1]}{\prod^{\bullet}}_{i=0}^{\lfloor log_2n\rfloor-1}y^{2^i}$.
For $2\nmid n$, it is easy to see that $|W|=2n+\lfloor log_2n\rfloor-1$ and $W$ is $n\mathbb{N}$-$1$-product free.
Thus it suffices to prove that $s_{n}(D_{2n})\leq 2n+\lfloor log_2n\rfloor$.

Fix a sequence $S$ of length $2n+\lfloor log_2n\rfloor$.
Let $\varphi_{\alpha}$ be a homomorphism mapping: $D_{2n}\rightarrow D_{2n}$ with $\varphi_{\alpha}(x)=xy^{\alpha}$ and $\varphi_{\alpha}(y)=y$, where $\alpha$ is an integer.
It is easy to see that for any integer $\alpha$, $T$ is a $1$-product subsequence of $S$ iff
$\varphi_{\alpha}(T)$ is a $1$-product subsequence of $\varphi_{\alpha}(S)$.
It follows that if $\varphi_{\alpha}(S)$ has a $1$-product subsequence $\varphi_{\alpha}(T)$ with $|\varphi_{\alpha}(T)|\equiv 0$ ($mol$ $n$),
then $T$ is a $1$-product subsequence of $S$ with $|T|\equiv 0$ ($mol$ $n$).
Hence, we can suppose that the term with the maximal multiplicity in $S_N$ is $x$.

Set
\begin{gather}
S_N=x^{[h(S_N)]}(a_1,a_1), \ldots ,(a_r,a_r),T_1,(c_1, \ldots ,c_u),  \label{S_N}
\end{gather}
where $h(S_N)$ is the maximal multiplicity in $S_N$ and
$T':=T_1(c_1, \ldots ,c_u)$ is squarefree with a maximal $1$-product subsequence $T_1$ ($T_1$ may be empty) and
a $1$-product free subsequence $T'(T_1)^{-1}$.
Obviously, $|T'|\leq n-1$ (Since $x\nmid T'$), $|T_1|$ is even and $c_1,\ldots, c_u$ are distinct.
By Corollary \ref{dihedral}, we have $u\leq 2\lfloor log_2n\rfloor+1$.
In addition, set
\begin{gather}
S_H=(b_1,b_1), \ldots ,(b_s,b_s),T_2,(d_1, \ldots ,d_v),  \label{S_H}
\end{gather}
where $T'':=T_2(d_1, \ldots ,d_v)$ is squarefree and
$T_2$ is a maximal $\pm$-$1$-product subsequence with $2||T_2|$ in $T''$ ($T_2$ may be empty).
Obviously, $d_1,\ldots, d_v$ are distinct.
By Corollary \ref{dihedral}, we have $v\leq 2\lfloor log_2n\rfloor+1$.
Let
\begin{gather}
S_0=(c_1, \ldots ,c_u)(d_1, \ldots ,d_v)=(xy^{\alpha_1}, \ldots ,xy^{\alpha_u})(y^{\beta_1}, \ldots ,y^{\beta_v}).  \label{S_0}
\end{gather}
Thus $|S_0|=u+v\leq 4\lfloor log_2n\rfloor+2$, $|S(S_0)^{-1}|\geq 2n-3\lfloor log_2n\rfloor-2$ and $|T_1|+|T_2|+|S_0|\leq 2n-1$ (Since $x\nmid T_1S_0$).

If there exists a $1$-product subsequence of length $n$ in $S_H$, then we are done.
Hence, we can suppose that $S_H$ is $n$-$1$-product free.
By Lemma \ref{basic} $(b)$, we must have that $|S_H|\leq 2n-3$ and then $|S_N|\geq \lfloor log_2n\rfloor+2\geq 3$.
If $h(S_N)=1$, then $S_N$ is squarefree and suppose $S_N=(c'_1,c'_2, \ldots ,c'_{u'})$.
It follows that $c'_1c'_2, \ldots ,c'_1c'_{u'}$ are distinct.
Thus $|\sum_2(c'_1,c'_2, \ldots ,c'_{u'})|\geq u'-1$.
Since $|S_H|=|S|-|S_N|=2n+\lfloor log_2n\rfloor-u'$ and $S_H$ is $n$-$1$-product free,
Lemma \ref{lem12} implies that $|\sum_{n-2}(S_H)|\geq n+\lfloor log_2n\rfloor-u'+1$.
It follows that
$$
|{\sum}_{n-2}(S_H)|+|{\sum}_{2}(c'_1,c'_2, \ldots ,c'_{u'})|\geq n+\lfloor log_2n\rfloor > n.
$$
By Lemma \ref{lem11} we have
$$
1\in {\sum}_{n-2}(S_H)+{\sum}_{2}(c'_1,c'_2, \ldots ,c'_{u'})\subseteq{\sum}_n (S).
$$

Suppose now that $h(S_N)\geq 2$.
If $S_H$ has a minimal $\pm$-$1$-product subsequence $T$ of odd length,
then by Lemma \ref{ODD} one can obtain that $|T|\leq n$.
If $|T|=n$, then Lemma \ref{ODD} implies that $T=g^{[n]}$, where $g$ is a generator of $H$, that is,
$T$ is a $1$-product subsequence of length $n$.
Suppose that $|T|\leq n-2$.
It will not be too confusing to assume that
\begin{gather*}
(ST^{-1})_N=S_N=x^{[h(S_N)]},(a_1,a_1), \ldots ,(a_r,a_r),T_1,(c_1, \ldots ,c_u),
\end{gather*}
with $h(S_N)\geq 2$,
and
\begin{gather*}
(ST^{-1})_H=(b_1,b_1), \ldots ,(b_s,b_s),T_2,(d_1, \ldots ,d_v).
\end{gather*}
with symbols having the same meaning as (\ref{S_H}).
Thus we can suppose that $|T|+h(S_N)+2r+2s\leq n-1$, since otherwise there exist $2h\in [1,h(S_N)]$, $r_1\in [0,r]$ and $s_1\in [0,s]$ such that
$$x^{[2h]},(a_1,a_1), \ldots ,(a_{r_1},a_{r_1}),(b_1,b_1), \ldots ,(b_{s_1},b_{s_1}),T$$
is a $1$-product subsequence of length $n$.
Hence, $|T_1|+|T_2|+u+v=|S|-(|T|+h(S_N)+2r+2s)\geq n+1+\lfloor log_2n\rfloor$.
Let $$A=\{i\in [0,n-1]: xy^i|T_1(c_1, \ldots ,c_{u})\} \ and\ B=\{j\in [0,n-1]: y^j|T_2(d_1, \ldots ,d_v)\}.$$
Thus $|A\cup B|\leq n$.
Since $T_1(c_1, \ldots ,c_{u})$ and $T_2(d_1, \ldots ,d_v)$ are squarefree,
we have that $n\geq |A\cap B|=|A|+|B|-|A\cup B|=|T_1|+u+|T_2|+v-|A\cup B|\geq 1+\lfloor log_2n\rfloor$.
Combining the above results yields that $|T|+h(S_N)+2r+2s+2|A\cap B|=|S|+|A\cap B|-|A\cup B|\geq n+1+2\lfloor log_2n\rfloor$.
Set $S_{AB}=\prod^{\bullet}_{i\in A\cap B}(y^i,xy^i)$.
Thus we can easily find a $1$-product subsequence of length $n$ in
$$x^{[h(S_N)]},(a_1,a_1), \ldots ,(a_{r},a_{r}),(b_1,b_1), \ldots ,(b_{s},b_{s}),T,S_{AB}.$$

In the following, let $h(S_N)\geq 2$ and suppose that $S_H$ does not have $\pm$-$1$-product subsequences of odd length.
Thus by Lemma \ref{ODD} we must have that $|S_H|\leq n-1$ and $|S_N|=|S|-|S_H|\geq n+1+\lfloor log_2n\rfloor$.
It follows from (\ref{S_H}) that
$(d_1, \ldots ,d_v)$ does not have $\pm$-$1$-product subsequences of even length.
Combining this with the assumption of $S_H$ yields that
\begin{gather}
(d_1, \ldots ,d_v) \label{d_1v}
\end{gather}
is $\pm$-$1$-product free.
This continues to imply that $v\leq \lfloor log_2n\rfloor$ by Lemma \ref{cyweight}.
Again by the assumption of $S_H$, we have that $1\nmid S_H$ which implies that $|T_2|+v\leq n-1$.
For $|T_1|+u\leq n-1$, we have $|T_1|+u+|T_2|+v\leq 2n-2$.

Firstly, suppose that $|T_1|+u+|T_2|+v=2n-2$.
Set $$A=\{i\in [1,n-1]: xy^i|T_1(c_1, \ldots ,c_{u})\}\ and\ B=\{j\in [1,n-1]: y^j|T_2(d_1, \ldots ,d_v)\}.$$
Since $x\nmid T_1(c_1, \ldots ,c_{u})$ and $1\nmid S_H$, we have $|A\cup B|\leq n-1$.
It follows that $n-1\geq |A\cap B|=|A|+|B|-|A\cup B|=|T_1|+u+|T_2|+v-|A\cup B|\geq n-1$,
i.e., $|A\cap B|=n-1$.
It is easy to see that
$$(x,x){\prod_{i\in A\cap B}}^{\bullet}(y^i,xy^i)$$
is a $1$-product subsequence of length $2n$.

Secondly, suppose that $|S_0|\leq \lfloor log_2n\rfloor$.
Then $|S_N(S_0)^{-1}|\geq n+1$ and $|S(S_0)^{-1}|\geq 2n$.
Set $\overline{S}=S(S_0)^{-1}$, if $|S(S_0)^{-1}|$ is even;
set $\overline{S}=S(S_0,x)^{-1}$, if $|S(S_0)^{-1}|$ is odd.
For the latter case, we must have that $|S(S_0)^{-1}|\geq 2n+1$.
Hence, $|\overline{S}|$ is even and $|\overline{S}|\geq 2n$.
By (\ref{S_N}) and (\ref{S_H}), one can easily obtain that
$\overline{S}_N$ is $1$-product and $\overline{S}_H$ is $\pm$-$1$-product
satisfying that $|\overline{S}_N|\geq 2$ and $|\overline{S}_N|$, $|\overline{S}_H|$ are even.
Combining them with Lemma \ref{dihedral-length} yields that $\overline{S}$ is $1$-product.
If $|\overline{S}((b_1,b_1), \ldots ,(b_s,b_s))^{-1}|\leq 2n$,
then there exists $s_1\in [0,s]$ such that
$$\overline{S}((b_1,b_1), \ldots ,(b_{s_1},b_{s_1}))^{-1}$$
is a $1$-product subsequence of length $2n$.
Suppose that $|\overline{S}((b_1,b_1), \ldots ,(b_s,b_s))^{-1}|> 2n$.
Since $|T_1|+|T_2|\leq |T_1|+u+|T_2|+v\leq 2n-2$,
there exists $r_1\in [0,r]$ and $\ell\in [1,h(S_N)]$ such that
$$x^{[\ell]},(a_1,a_1), \ldots ,(a_{r_l},a_{r_l}),T_1,T_2$$
is a $1$-product subsequence of length $2n$.

Finally, by the above analysis, it suffices to prove that
if $h(S_N)\geq 2$, $S_H$ does not have $\pm$-$1$-product subsequences of odd length, $|T_1|+u+|T_2|+v\leq 2n-3$ and $|S_0|\geq \lfloor log_2n\rfloor+1$,
then there exists a $1$-product subsequence of length $2n$ in $S$.
In the following, we distinguish two cases in terms of whether or not $h(S_N)\geq \lfloor log_2n\rfloor+1$.

\textit{Case $1$:} $h(S_N)\geq \lfloor log_2n\rfloor+1$.

Let $\psi$ be a map $D_{2n}\rightarrow C_n$ with $\psi(xy^{\alpha})=\alpha$ and $\psi(y^{\beta})=\beta$.
It follows from (\ref{S_0}) that
$$\psi(S_0)=(\alpha_1, \ldots ,\alpha_u)(\beta_1, \ldots ,\beta_v)$$
with $|\psi(S_0)|=|S_0|\geq \lfloor log_2n\rfloor+1$.
By Lemma \ref{cyweight}, we have that
$$\psi(S_0)=\psi(L_1), \ldots ,\psi(L_k),\psi(L')$$
where $k\geq 1$, $\psi(L_i)$ is minimal $\pm$-zero-sum
with $|L_i|=|\psi(L_i)|\leq \lfloor log_2n\rfloor+1$ for $1\leq i\leq k$
and $|L'|=|\psi(L')|\leq \lfloor log_2n\rfloor$.
Suppose that $L_i=L'_iL''_i$ with $\sigma_{\pm}(\psi(L_i))=\sigma(\psi(L'_i))-\sigma(\psi(L''_i))=0$.
In addition, we have $L'_i=L'_{iH}L'_{iN}$ and $L''_i=L''_{iH}L''_{iN}$,
where $L'_{iH}L''_{iH}=L_i\cap H=L_{iH}$ and $L'_{iN}L''_{iN}=L_i\cap N=L_{iN}$.
By a reasonable rearrangement for the signals $\pm 1$ appearing in $\psi(L'_i)$ and $\psi(L''_i)$,
we have that for any $k_1\in [1,k]$
$$|\sum_{i=1}^{k_1}|L'_{iN}|-\sum_{i=1}^{k_1}|L''_{iN}||\leq \lfloor log_2n\rfloor+1$$
and $\Delta:=\sum_{i=1}^{k}|L'_{iN}|-\sum_{i=1}^{k}|L''_{iN}|\geq 0$.
Since $h(S_N)\geq \lfloor log_2n\rfloor+1\geq \Delta$,
we have
$$S'_{0N}:=L'_{1N}, \ldots ,L'_{kN}$$
and
$$S''_{0N}:=x^{[\Delta]},L''_{1N}, \ldots ,L''_{kN}$$
with $|S'_{0N}|=|S''_{0N}|$.
Since $\psi(L_i)$ is minimal $\pm$-zero-sum for $1\leq i\leq k$,
it is easy to see that
\begin{align*}
&\prod_{xy^{\alpha_i}|S'_{0N},xy^{\alpha_j}|S''_{0N}}
(xy^{\alpha_j}\cdot xy^{\alpha_i})\prod_{i=1}^{k}\pi(L'_{iH})\prod_{i=1}^{k}(\pi(L''_{iH}))^{-1} \\
&=y^{\sigma(\psi(S'_{0N}))-\sigma(\psi(S''_{0N}))}\cdot y^{\sum_{i=1}^{k}(\sigma(\psi(L'_{iH}))-\sigma(\psi(L''_{iH})))}\\
&=y^{\sum_{i=1}^{k}(\sigma(\psi(L'_{iH}))+\sigma(\psi(L'_{iN}))-\sigma(\psi(L''_{iH}))-\sigma(\psi(L''_{iN})))}\\
&=y^{\sum_{i=1}^{k}(\sigma(\psi(L'_{i}))-\sigma(\psi(L''_{i})))}
=y^{\sum_{i=1}^{k}(\sigma_{\pm}(\psi(L_i)))}=1.
\end{align*}
It follows that
\begin{gather}
\overline{S}:=x^{[\Delta]},L_1, \ldots ,L_k,T_1,T_2  \label{overlineS}
\end{gather}
is $\pm$-$1$-product with $(L_1, \ldots ,L_k,T_1,T_2)$ squarefree.
If $|\overline{S}|=2n$, then $|T_1|+u+|T_2|+v\leq 2n-3$ implies that $|\overline{S}_N|\geq \Delta>0$.
Combining this with Lemma \ref{dihedral-length} yields that $\overline{S}$ is a $1$-product subsequence of length $2n$,
we are done.
Suppose that $|\overline{S}|<2n$.
For $|L'|\leq \lfloor log_2n\rfloor$, we have $|S(L')^{-1}|=|\overline{S}|+h(S_N)-\Delta+2r+2s\geq 2n$.
Since $|(S(L')^{-1})_N|>0$, one can obtain that $S(L')^{-1}$ is $1$-product.
Thus there exists $\triangle'\in[0,h(S_N)-\Delta]$, $r_1\in [0,r]$ and $s_1\in [0,s]$ such that
$$\overline{S}':=\overline{S},x^{[\triangle']},(a_1,a_1), \ldots ,(a_{r_1},a_{r_1}),(b_1,b_1), \ldots ,(b_{s_1},b_{s_1})$$
satisfies that $|\overline{S}'|=2n$ and $|(\overline{S}')_N|>0$.
It follows that $\overline{S}'$ is a $1$-product subsequence of length $2n$.

Now suppose $|\overline{S}|>2n$.
Since $|T_1|+|T_2|+|S_0|=|T_1|+u+|T_2|+v\leq 2n-3$,
by (\ref{overlineS}) we must have $\triangle\geq 4$.
Combining this with $\Delta\leq \lfloor log_2n\rfloor+1$ yields that $n\geq 9$.
In (\ref{overlineS}) we suppose that $$\overline{S}_1:=x^{[\Delta]},L_1, \ldots ,L_k.$$
Combining the above results with (\ref{S_N}), (\ref{S_H}) yields that $\overline{S}_1$ is $\pm$-$1$-product with
$|\overline{S}_1|\leq \Delta+u+v\leq 5\lfloor log_2n\rfloor +3\leq 2n$.
In addition, it is easy to see that $\overline{S}$ and $\overline{S}_1$ are $1$-product.
Set
\begin{equation}
\label{overlineS''}
   \begin{aligned}
   &\overline{S}''=\overline{S},x^{[h(S_N)-\Delta-\sigma]},(a_1,a_1), \ldots ,(a_{r},a_{r}),(b_1,b_1), \ldots ,(b_{s},b_{s})  \\
   &\overline{S}''_1=\overline{S}_1,x^{[h(S_N)-\Delta-\sigma]},(a_1,a_1), \ldots ,(a_{r},a_{r}),(b_1,b_1), \ldots ,(b_{s},b_{s}) \\
   \end{aligned}
  \end{equation}
where $\sigma=0$ if $h(S_N)-\Delta$ is even; $\sigma=1$ if $h(S_N)-\Delta$ is odd.
Thus $\overline{S}''$ and $\overline{S}''_1$ are also $1$-product with $\overline{S}''=\overline{S}''_1T_1T_2$ and $|\overline{S}''|\geq |\overline{S}|>2n$.
In addition, we can suppose that $|\overline{S}''_1|< 2n$, since otherwise by $|\overline{S}_1|\leq 2n$
we can easily find a $1$-product subsequence of length $2n$ in $\overline{S}''_1$.

Since $|\overline{S}|>2n$ and $\Delta\leq \lfloor log_2n\rfloor+1$,
from (\ref{S_0}) and (\ref{overlineS}) one can obtain that
$\overline{S}(x^{[\Delta]})^{-1}=T_1T_2S_0(L')^{-1}$ is squarefree with
$|\overline{S}(x^{[\Delta]})^{-1}|=|T_1|+|T_2|+u+v-|L'|\geq 2n-\lfloor log_2n\rfloor$.
Since $|T_1|+u\leq n-1$ and $|T_2|+v \leq n-1$,
we have that $n-\lfloor log_2n\rfloor+1\leq|T_1|+u\leq n-1$
and $n-\lfloor log_2n\rfloor+1\leq|T_2|+v\leq n-1$.
Suppose that
$$
T_1,(c_1, \ldots ,c_{u})=C_1, \ldots ,C_{u_1},
$$
and
$$
T_2,(d_1, \ldots ,d_v)=D_1, \ldots ,D_{v_1},
$$
where the power exponents of $y$ in each $C_i$, $D_j$ are continuous and
$$max\{t: xy^t|C_i\}+1<min\{t':xy^{t'}|C_{i+1}\},$$
$$max\{l: y^l|D_j\}+1<min\{l':y^{l'}|D_{j+1}\}.$$
From this we can obtain two maximal subsequences:
$$F_1, \ldots ,F_{2\omega}|T_1(c_1, \ldots ,c_{u})$$
and
$$(G_1,G'_1), \ldots ,(G_\lambda,G'_\lambda)|T_2(d_1, \ldots ,d_v)$$
satisfying that $|F_i|=2$ and $\pi(F_i)=y$ for $1\leq i\leq 2\omega$;
$G_j=\{y^{\alpha_j},y^{\alpha_{j}+1}\}$ and $G'_j=\{y^{\beta_j},y^{\beta_{j}+1}\}$ for $1\leq j\leq \lambda$.
It is easy to see
$$n-\lfloor log_2n\rfloor+1\leq |T_1|+u=\sum_{i=1}^{u_1}|C_i|\leq n-u_1+1$$
and
$$n-\lfloor log_2n\rfloor +1\leq |T_2|+v=\sum_{j=1}^{v_1}|D_j|\leq n-v_1+1$$
Thus
$$4\omega=\sum_{i=1}^{2\omega}|F_i|\geq \sum_{i=1}^{u_1}2(\lfloor\frac{|C_i|}{2}\rfloor)-2\geq \sum_{i=1}^{u_1}(|C_i|-1)-2=|T_1|+u-u_1-2\geq n-2\lfloor log_2n\rfloor-1.$$
Similarly,
$$4\lambda\geq \sum_{j=1}^{v_1}(|D_j|-1)-2=|T_2|+v-v_1-2\geq n-2\lfloor log_2n\rfloor-1.$$
In addition, ${\prod}_{i=1}^{\omega}(\pi(F_{2i-1})\cdot \pi(F_{2i}^{-1}))=1$
and ${\prod}_{j=1}^{\lambda}((y^{\alpha_j}\cdot (y^{\alpha_{j}+1})^{-1})\cdot ((y^{\beta_j})^{-1}\cdot y^{\beta_{j}+1}))=1$.
Thus $F_1, \ldots ,F_{2\omega}$ is $1$-product and
$(G_1,G'_1), \ldots ,(G_\lambda,G'_\lambda)$ is $\pm$-$1$-product.
From (\ref{S_N}) and (\ref{d_1v}) we know that $(c_1, \ldots ,c_{u})$ is $1$-product free
and $(d_1, \ldots ,d_v)$ is $\pm$-$1$-product free.
Hence, we can suppose that
$$F_1, \ldots ,F_{2\omega}|T_1$$
and
$$(G_1,G'_1), \ldots ,(G_\lambda,G'_\lambda)|T_2.$$

In (\ref{overlineS''}), suppose $h(S_N)-\Delta-\sigma+2r+2s\geq 2$.
It is easy to see that
$$\overline{S}''':=\overline{S}''(F_1, \ldots ,F_{2\omega},(G_1,G'_1), \ldots ,(G_\lambda,G'_\lambda),(z,z))^{-1}$$
is a $1$-product subsequence where $(z,z)|\overline{S}''(\overline{S})^{-1}$.
In addition,
$|\overline{S}'''|=|\overline{S}''|-4\omega-4\lambda-2\leq |S|-4\omega-4\lambda-2\leq 5\lfloor log_2n\rfloor\leq 2n$.
For $|\overline{S}''|>2n$,
there is $\sigma=0$ or $1$, $\omega_1\in [0,\omega]$ and $\lambda_1\in [0,\lambda]$ such that
$|\overline{S}'''|+4\omega_1+4\lambda_1+2\sigma=2n$
where $\sigma=0$, if $2n-|\overline{S}'''|\equiv 0$ (mod $4$); $\sigma=1$, if $2n-|\overline{S}'''|\equiv 2$ (mod $4$).
Thus,
$$\overline{S}''',F_1, \ldots ,F_{2\omega_1},(G_1,G'_1), \ldots ,(G_{\lambda_1},G'_{\lambda_1}),(z,z)^{\sigma}$$
is a $1$-product subsequence of length $2n$.

If $h(S_N)-\Delta-\sigma+2r+2s\leq 1$,
then $h(S_N)-\Delta-\sigma\leq 1$ and $r=s=0$.
It follows that $$S=x^{[h(S_N)]},T_1(c_1, \ldots ,c_{u})T_2(d_1, \ldots ,d_v).$$
From (\ref{overlineS''}), we know that $h(S_N)-\Delta-\sigma$ is even.
Thus $h(S_N)-\Delta-\sigma=0$.
For $\Delta\leq \lfloor log_2n\rfloor+1$, we have
$h(S_N)=\Delta+(h(S_N)-\Delta-\sigma)+\sigma\leq\lfloor log_2n\rfloor+2$.
Thus $|T_1|+u+|T_2|+v=|S|-h(S_N)\geq 2n-2$.
This is in contradiction to $|T_1|+u+|T_2|+v\leq 2n-3$.

\textit{Case 2:} $h(S_N)\leq \lfloor log_2n\rfloor$.

By (\ref{S_N}),
we have that $\psi(S_N)=0^{[h(S_N)]}\psi(S_{N1})$,
where $\psi$ is a map: $D_{2n}\rightarrow C_n$ with $\psi(xy^{\alpha})=\alpha$, $\psi(y^{\beta})=\beta$,
and
$S_{N1}=S_N(x^{[h(S_N)]})^{-1}$.
For $|S_N|\geq n+\lfloor log_2n\rfloor+1$,
we have that $|S_{N1}|=|\psi(S_{N1})|\geq n+\lfloor log_2n\rfloor+1-h(S_N)\geq n$.
By using Lemma \ref{maxlength} on $T$ repeatedly,
there exist minimal $1$-product subsequences $\psi(\overline{T_1}), \ldots ,\psi(\overline{T_{k'}})$ in $\psi(S_{N1})$
with $|\psi(\overline{T_i})|\leq h(S_N)$ for $1\leq i\leq k'$
and $\lfloor log_2n\rfloor-h(S_N)+1\leq|\psi(\overline{T_1})|+ \ldots +|\psi(\overline{T_{k'}})|\leq \lfloor log_2n\rfloor$.
It will not be too confusing to suppose that
\begin{gather*}
S_N=x^{[h(S_N)]},\overline{T_1}, \ldots ,\overline{T_{k'}},(a_1,a_1), \ldots ,(a_r,a_r),T_1,(c_1, \ldots ,c_u)
\end{gather*}
with $h(S_N)+|\overline{T_1}|+ \ldots +|\overline{T_{k'}}|\geq \lfloor log_2n\rfloor+1$,
and
\begin{gather*}
S_H=(b_1,b_1), \ldots ,(b_s,b_s),T_2,(d_1, \ldots ,d_v).
\end{gather*}
In addition, we can suppose that
\begin{gather*}
S_0=(c_1, \ldots ,c_u)(d_1, \ldots ,d_v)
\end{gather*}
with $|S_0|\geq \lfloor log_2n\rfloor+1$
and
$$\psi(S_0)=\psi(L_1), \ldots ,\psi(L_k),\psi(L').$$
Note that the symbols have the same meaning with the above.
By imitating the proof of Case $1$, we have that $L_{iN}=L'_{iN}L''_{iN}$
and for any $k_1\in [1,k]$
$$|\sum_{i=1}^{k_1}|L'_{iN}|-\sum_{i=1}^{k_1}|L''_{iN}||\leq \lfloor log_2n\rfloor+1$$
and $\Delta:=\sum_{i=1}^{k}|L'_{iN}|-\sum_{i=1}^{k}|L''_{iN}|\geq 0$.

If $\Delta\geq |\overline{T_1}|+ \ldots +|\overline{T_{k'}}|$,
then there is $\Delta_1\in [0,\Delta]$ such that
$$L_1, \ldots ,L_k,\overline{T_1}, \ldots ,\overline{T_{k'}},(\underbrace{x, \ldots ,x} \limits_{\Delta_1})$$
is $1$-product, where
$$\sum_{i=1}^{k}|L'_{iN}|=\sum_{i=1}^{k}|L''_{iN}|+|\overline{T_1}|+ \ldots +|\overline{T_{k'}}|+\Delta_1.$$

If $\Delta< |\overline{T_1}|+ \ldots +|\overline{T_{k'}}|$,
then there are $\overline{T'}=\overline{T_1}, \ldots ,\overline{T_{k'_1}}$
and $\overline{T''}=\overline{T_{k'_1}}, \ldots ,\overline{T_{k'}}$
such that
$$|(\sum_{i=1}^{k}|L'_{iN}|+|\overline{T'}|)-(\sum_{i=1}^{k}|L''_{iN}|+|\overline{T''}|)|\leq h(S_N).$$
It follows that there is $\Delta_2\in [0,h(S_N)]$ such that
$$L_1, \ldots ,L_k,\overline{T_1}, \ldots ,\overline{T_{k'}},(\underbrace{x, \ldots ,x} \limits_{\Delta_2})$$
is $1$-product.

Repeat the reasoning in Case $1$ and
we can find a $1$-product subsequence of length $2n$.
The proof is complete.

\qed

\textit{The proof of Theorem \ref{main}:}
Let $D_{2n}=\langle x,y:x^2=y^n=1, xy=y^{-1}x\rangle$.

Let $d$ be odd and $n|d$. Suppose that $d=kn$ for some positive integer $k$ and then $2\nmid n$, $2\nmid k$.
Set $W=x^{[2d-1]}{\prod^{\bullet}}_{i=0}^{\lfloor log_2n\rfloor-1}y^{2^i}$.
It is easy to see that $W$ is $d\mathbb{N}$-$1$-product free and $|W|=2d+\lfloor log_2n\rfloor-1$.
Thus it suffices to prove that $s_{d\mathbb{N}}(D_{2n})\leq 2d+\lfloor log_2n\rfloor$.

Fix a sequence $S$ of length $2d+\lfloor log_2n\rfloor$ in $D_{2n}$.
By using Lemma \ref{basic} $(g)$ on $S$ repeatedly,
we have a decomposition:
$$S=T_1\ldots T_{k-1}S_1$$
where $|S_1|=2n+\lfloor log_2n\rfloor$
and each $T_i$ is a $1$-product subsequence of length $2n$.
Theorem \ref{n-dihedral} implies that $S_1$ has a $1$-product subsequence $S_2$ of length $n$ or $2n$.
If $|S_2|=2n$, then $T_1\ldots T_{k-1}S_2$ is a $1$-product subsequence of length $2d$.
If $|S_2|=n$, then $T_1\ldots T_{\frac{k-1}{2}}S_2$ is a $1$-product subsequence of length $d$.
The proof of the first assertion is complete.

Let $gcd(n,d)=1$.
It is easy to see that $W=xy^{[nd-1]}$ is a $d\mathbb{N}$-$1$-product sequence of length $nd$.
Thus it suffices to prove that $s_{d\mathbb{N}}(D_{2n})\geq nd+1$.

Fix a sequence $S$ of length $nd+1$ in $D_{2n}$.
By using Lemma \ref{lem40} $(ii)$ on $S$ repeatedly,
we have a decomposition:
$$S=T_1 \ldots T_{d}T$$
where $\pi (T_i)=1$, $1\leq |T_i|\leq n$ for any $i\in [1,d]$, and $|T|\geq 1$.
Let $$T=|T_1|, \ldots ,|T_d|$$
be a sequence of length $d$ in the cyclic group $C_d$.
By Lemma \ref{basic} $(a)$, $T$ has a zero-sum subsequence
$$T'=|T_{i_1}|, \ldots ,|T_{i_k}|,$$
that is, $\sum_{j=1}^{k}|T_{i_j}|\equiv 0$ (mod $d$).
Hence, $$S'=T_{i_1}, \ldots ,T_{i_k}$$
is a $1$-product subsequence
with length $|S'|\equiv 0$ (mod $d$) and we complete the proof.

\qed

\section{The proof of Theorem \ref{thm2}}

Let $G_{pq}=\langle x,y: x^p=y^q = 1,yx = xy^s,ord_q(s) =p, \ and \ p,q\ primes\rangle$ be an nonabelian group.
Thus we must have that $p\geq 2$ and $p|q-1$.
Note that for $p=2$, this is simply the dihedral group of order $2q$.
Therefore, in the following we assume $p\geq 3$, so $q\geq 2p+1$.

Similarly to the previous sections, let $H$ be the cyclic subgroup of $G_{pq}$ generated by $y$,
and let $N=G_{pq}\setminus H$.
We can further split up $N$ into cosets of $H$: let $N_i=x^iH$, so $N=N_1\cup \ldots\cup N_{p-1}$.
Note that these cosets (with $H$ included) form a group isomorphic to $\mathbb{Z}_p$,
and since $D(\mathbb{Z}_p)=p$, any sequence of $p$ elements in $G_{pq}$ must have a nonempty subsequence $T$ with $\pi(T)\in H$.
If $A$ and $B$ are subsets of a group, then we define $A+B=\{ab:a\in A, b\in B\}$.
In the preliminaries, we define $\Pi(S)=\{\prod_{i\in [1,\ell]}g_{\tau(i)} : \tau \ is\ a\ permutation\ of\ [1,\ell]\}$
for a sequence $S=g_1 \cdot\ldots\cdot g_{\ell}\in \mathscr{F}(G)$.

The following two lemmas will be repeatedly used in the proof of our main theorem.

\begin{lem}[Cauchy-Davenport inequality \cite{[N]}, pp 44-45] \label{Cauchy-Davenport}
For $q$ prime, and for any $k$ nonempty sets $A_1,\ldots, A_k\subseteq \mathbb{Z}_q$,
$$|A_1+\ldots+A_k|\geq min\{q,|A_1|+\ldots+|A_k|-k+1\}.$$
In particular, if $A$ and $B$ are two nonempty subsets of $\mathbb{Z}_q$, then $|A+B|\geq min \{q,|A|+|B|-1\}$.
\end{lem}

\begin{lem}[\cite{[B]}, Theorem 15] \label{lem1}

If $p\geq 3$ and $p|q-1$, then $\mathsf{s}_{pq}(G_{pq})=pq+p+q-2$.
\end{lem}

By imitating the proof of Theorem 15 of \cite{[B]}, we can prove Theorem \ref{thm2}.
In the proof, we will cite some results showed in Theorem 15 of \cite{[B]} and ignore their proofs.

\textit{The proof of Theorem \ref{thm2}:}
Let $G_{pq}=\langle x,y: x^p=y^q = 1,yx = xy^s,ord_q(s) =p, \ and \ p,q\ primes\rangle$,
$H=\langle y\rangle$, and $N=G_{pq}\setminus H$.
Put $d=lcm[kp,q]+gcd(kp,q)-1$.
It is easy to see that $W=x^{[p-1]}y^{[d-1]}$ is a $kp\mathbb{N}$-$1$-product sequence of length $p+d-2$.
Thus it suffices to prove that $s_{kp\mathbb{N}}(G_{pq})\geq p+d-1$.

Fix a sequence $S$ of length $p+d-1$ in $G_{pq}$.
If $|S\cap N|\leq p-1$, then $|S\cap H|\geq d$.
Combining Lemma \ref{basic} $(c)$ and $(d)$ yields that $D(H\oplus C_{kp})=D(C_q\oplus C_{kp})=D(C_{gcd(kp,q)}\oplus C_{lcm(kp,q)})=lcm(kp,q)+gcd(kp,q)-1=d=D^*(H\oplus C_{kp})$.
By Lemma \ref{basic} $(e)$, we have that $\mathsf{s}_{kp\mathbb{N}}(H)=d$ which implies that $S\cap H$ has a $1$-product subsequence $T$ satisfying that $|T|\equiv 0$ (mod $kp$).
If $q|k$, then $|S|=kp+p+q-2$.
By using Lemma \ref{lem1} on $S$ repeatedly,
we have a decomposition:
$$S=T_1\ldots T_{\frac{k}{q}}S_1$$
where $|S_1|=p+q-2$
and each $T_i$ is a $1$-product subsequence of length $pq$.
It follows that $(T_1, \ldots ,T_{\frac{k}{q}})$ is a $1$-product subsequence of length $kp$.

Now suppose that $|S\cap N|\geq p$ and $(q,k)=1$.
Thus $|S|=kpq+p-1$.

In the following, by imitating the proof of Theorem 15 of \cite{[B]}, we can complete the proof.

In Theorem 15 of \cite{[B]}, Bass showed the following three claims:

\textbf{Claim $1$:}

For any sequence $S$ of length $p+d-1$ in $G_{pq}$, we have the following decomposition:
$$
S=(A_1, \ldots ,A_r)(A'_1, \ldots ,A'_{r'})(F_1, \ldots ,F_{\ell})(F'_1, \ldots ,F'_{\ell'})E
$$
where
\begin{itemize}
  \item For all $1\leq k\leq r$, $A_k$ is a subsequence of $S\cap N_i$ for some $i$ and has $\pi(A_k)=1$, $|A_k|=p$.
  For all $1\leq k\leq r'$, $A'_k$ has $\pi(A'_k)=1$ and $|A'_k|=p$.
  \item For all $1\leq k\leq \ell$, $F_k$ is a subsequence of $S\cap N_i$ for some $i$, $\pi(F_k)=y^{m_k}\neq 1$, $|F_k|=p$, and $\Pi(F_k)=\{y^{m_k},y^{m_ks},\ldots ,y^{m_ks^{p-1}}\}$.
      For all $1\leq k\leq \ell'$, $F'_k$ has $\pi(F'_k)=y^{m'_k}\neq 1$ and $|F'_k|=p$.
  \item $|(S((A_1, \ldots ,A_r)(F_1, \ldots ,F_{\ell}))^{-1})_N|\leq (p-1)^2$ and $|E|\leq 2p-2$.
  Note that $|E|\equiv |S|=kpq+p-1\equiv -1$ (mod $p$), since $p|q-1$.
  Since $0\leq |E|\leq 2p-2$, this implies $|E|=p-1$, and therefore $(\ell+\ell'+r+r')p=kpq$.
\end{itemize}

\textbf{Claim $2$:} (Lemma $16$ of \cite{[B]}) Let $Z=(z_1,\ldots,z_p)$ be a sequence of $p$ elements in $N_i$.
Then for any $k\in [1,\ell']$, the sequence $ZF'_1\ldots F'_{\ell'}$ attains all products in the set
$$
\{\Pi(Z)\}+\{y^{m'_1},y^{m'_1s}, \ldots ,y^{m'_1s^{p-1}}\}+ \ldots +\{y^{m'_k},y^{m'_ks}, \ldots ,y^{m'_ks^{p-1}}\}.
$$

\textbf{Claim $3$:} (Subcase $3A$ of \cite{[B]}) For any nonempty sequence $W\in \mathscr{F}(N)$ with $\pi(W)\in H$, a product of the sequence $P=(d_1,\ldots,d_{q-1})W$ can take on any value in $H$, where $d_i\in H$ for all $1\leq i\leq q-1$.

By Claim $1$, we can assume that $r+r'\leq k-1$, since otherwise there exists a $1$-product subsequence of length $kp$ in $(A_1, \ldots ,A_r)(A'_1, \ldots ,A'_{r'})$.
Thus, $(l+l')p=\sum_{i=1}^{\ell}|F_i|+\sum_{j=1}^{\ell'}|F'_j|=|S|-|E|-(\sum_{i=1}^{r}|A_i|+\sum_{j=1}^{r'}|A'_j|)=kpq-(r+r')p\geq kpq-(k-1)p=kp(q-1)+1\geq q$.

In the following, we distinguish three cases:

Case $1:$ $\ell\geq 1$.

Recall the definitions of $F_i$ and $\Pi(F_1)$ in Claim $1$.
Use Claim $2$ to $k=\ell'$, $Z=F_1$ and we have that the products of $F_1,F'_1, \ldots ,F'_{\ell'}$ attain all values in the set
$$
\{\Pi(F_1)\}+\{y^{m'_1},y^{m'_1s}, \ldots ,y^{m'_1s^{p-1}}\}+ \ldots +\{y^{m'_{\ell}},y^{m'_{\ell}s}, \ldots ,y^{m'_{\ell}s^{p-1}}\}
$$
where $\Pi(F_1)=\{y^{m_1},y^{m_1s}, \ldots ,y^{m_1s^{p-1}}\}$.
Since $\pi(A_i)=1$ and $\pi(A'_j)=1$ for all $1\leq i\leq r$ and $1\leq j\leq r'$,
the sequence $SE^{-1}=(A_1, \ldots ,A_r)$$(A'_1, \ldots ,A'_{r'})$ \\
$(F_1, \ldots ,F_{\ell})$$(F'_1, \ldots ,F'_{\ell'})$ can take on all products in a set of the form

\begin{align*}
X:&=\{y^{m_1},y^{m_1s}, \ldots ,y^{m_1s^{p-1}}\}+\ldots +\{y^{m_l},y^{m_ls}, \ldots ,y^{m_ls^{p-1}}\}   \\
&+\{y^{m'_1},y^{m'_1s}, \ldots ,y^{m'_1s^{p-1}}\}+ \ldots +\{y^{m'_l},y^{m'_ls}, \ldots ,y^{m'_ls^{p-1}}\}.
\end{align*}

By Lemma \ref{Cauchy-Davenport}, we have that $|X|\geq min\{q,(\ell+\ell')q\}=q$,
that is, $SE^{-1}$ is a $1$-product subsequence of length $kpq$.

Case $2:$ $\ell=0$ and $r\geq 1$.

By the definitions of $A_i$,
we can use Claim $2$ to $k=\ell'$ and $Z=A_1$.
From $\pi(A_1)=1$ it follows that the products of $A_1,F'_1, \ldots ,F'_{\ell'}$ attain all values in the set
$$
\{y^{m'_1},y^{m'_1s}, \ldots ,y^{m'_1s^{p-1}}\}+ \ldots +\{y^{m'_{\ell}},y^{m'_{\ell}s}, \ldots ,y^{m'_{\ell}s^{p-1}}\}.
$$
Repeat the reasoning of Case $1$ and one can obtain that $SE^{-1}$ is a $1$-product subsequence of length $kpq$.

Case $3:$ $l=r=0$.

This implies that $p\leq |S\cap N|\leq (p-1)^2$ and $|S\cap H|\geq kpq+p-1-(p-1)^2$.
Set $$S\cap H=1^{[n]}d_1,\ldots,d_v=1^{[n]}y^{m_1},\ldots,y^{m_v}$$
where $n=v_{1}(S\cap H)$ and $d_i=y^{m_i}$ is an nonidentity element for all $1\leq i\leq v$.
If $n\geq kp$, we are done.
If $n\leq kp-1$, then $v=|S\cap H|-n\geq kpq+p-(p-1)^2-kp\geq q-1$.
Consider $S\cap N$ as a sequence in $G_{pq}/H$ and suppose that
$W$ is its subsequence of maximum length satisfying $\pi(W)\in H$.
It follows that $|(S\cap N)W^{-1}|\leq p-1$.
For $|S\cap N|\geq p$, $W$ is nonempty.
By Claim $3$, we have that a product of the sequence $P=(d_1,\ldots,d_{q-1})W$ can take on any value in $H$.
Since $|P|=|W|+q-1\leq |S\cap N|+q-1\leq (p-1)^2+q-1\leq kpq$ and $|(S\cap H)W|=|S|-|(S\cap N)W^{-1}|\geq kpq$, take any subsequence $T$ of $(S\cap H)(d_1,\ldots,d_{q-1})^{-1}$ with length $|T|=kpq-|P|$ and let $\pi(P)=\pi(T)^{-1}$.
It follows that $TP$ is a $1$-product subsequence of length $kpq$ and the proof is complete.

\qed

\section*{Acknowledgments}
This work is supported by NSF of China  (Grant No. 11671153).
The author is sincerely grateful to the anonymous referee
for useful comments and suggestions.

\section*{References}


\begin{thebibliography}{99}

\bibitem{[AD]} S. Adhikari, Y. Chen, J. Friedlander, et al, Contributions to zero-sum problems, {\it Discrete Math.}, \textbf{306.1} (2006), 1-10.

\bibitem{[B]} J. Bass, Improving the Erd\H{o}s-Ginzburg-Ziv theorem for some non-abelian groups, {\it J. Number Theory}, \textbf{126.2} (2007), 217-236.

\bibitem{[GHPS]} W. Gao, D. Han, J. Peng, and F. Sun, On zero-sum subsequences of length
$k exp(G)$, {\it J. Combin. Theory Ser. A}, \textbf{125.1} (2014), 240-253.

\bibitem{[GHQQZ]} W. Gao, D. Han, G. Qian, Y. Qu, and H. Zhang, On additive bases
II, {\it Acta Arith.}, to appear.

\bibitem{[GL]} W. Gao and Z. Lu, The Erd\H{o}s-Ginzburg-Ziv theorem for dihedral groups, {\it J. Pure Appl. Algebra}, \textbf{212.2} (2008), 311-319.

\bibitem{[GGS]} A. Geroldinger, D. Grynkiewicz, and W. Schmid,
Zero-sum problems with congruence conditions, {\it Acta Math. Hungar.},
\textbf{131} (2011), 323-345.

\bibitem{[GH]} A. Geroldinger and F. Halter-Koch, {\it Non-unique factorizations,
Algebraic, Combinatorial and Analytic Theory}, Pure and Applied Mathematics, vol.278, Chapman and Hall/CRC, 2006.

\bibitem{[MBR]} F. E. B. Mart\'{\i}nez, S. Ribas, Extremal product-one free sequences in Dihedral and Dicyclic groups, {\it Discrete Math.}, \textbf{341.2} (2018), 570-578.

\bibitem{[N]} M. Nathanson, {\it Additive Number Theory: Inverse Problems and the Geometry of Sumsets}, Springer, 1996.

\bibitem{[QH]} Y. Qu and D. Han, An inverse theorem for additive bases, {\it Int. J. Number Theory}, (2015), 1-10.
\end{thebibliography}

\end{document}